\documentclass[11pt,letterpaper]{article}
\usepackage[hmargin=1in, vmargin=1in]{geometry}
\usepackage{graphicx}
\usepackage{amsmath}
\usepackage{dsfont}
\usepackage{indentfirst}
\usepackage{sidecap}
\usepackage{float}
\usepackage{multicol}
\usepackage{caption}
\usepackage{subcaption}
\usepackage[hyphens]{url}
\usepackage{hyperref}

\usepackage{fancyhdr}
\title{3D Printing an Octohedron}
\author{Edward Aboufadel \\ Grand Valley State University}
\date{Version 1.0, July 2014}

\begin{document}
\maketitle

\section{Introduction}

The purpose of this short paper is to describe a project to manufacture a regular octohedron on a 3D printer.  We assume that the reader is familiar with the basics of 3D printing.  In the project, we use fundamental ideas to calculate the vertices and faces of an octohedron.  Then, we utilize the OPENSCAD program to create a virtual 3D model and an STereoLithography (.stl) file that can be used by a 3D printer.  Finally, the object is created on a 3D printer.

\section{The Octohedron:  First Attempt}

An octohedron is a polyhedron with eight sides (all equilateral triangles) and six vertices.  It can be constructed in this way in ${\bf R}^3$:  start with a square of side length $s$ sitting in the $xy$-plane.  Identify the line normal to the plane that contains the center of the square.  Then, find the two points on the line (one above the plane and one below) so that the distance from each point to the four corners of the square equals $s$.  These two points and the four corners of the square are the six vertices of a regular octohedron, and each of the eight faces come from using one of the two points on the normal line as well as two adjacent corners of the square. See Figure 1.

In order to use OPENSCAD and to 3D print, we need to define specific coordinates for our points.  So, we start with a square with side length $s = 1$.  The four corners of the square we will use are:  $p_0 = (0.0, 0.0, 0.0)$, $p_1 = (1.0, 0.0, 0.0)$, $p_2 = (1.0, 1.0, 0.0)$ and $p_3 =(0.0, 1.0, 0.0)$.  (We start the indices at 0 because that is the standard when we use OPENSCAD.)

The normal line to this square is parallel to the $z$-axis through $(0.5, 0.5, 0)$, and for our two other vertices, we need to find points $(0.5, 0.5, \hat{z})$ that are a distance 1 from the four corners points of the square.  The distance formula yields $(0.5)^2 + (0.5)^2 + \hat{z}^2 = 1$, so $\hat{z}^2 = 0.5$, and $\hat{z} = \pm \sqrt{0.5} \approx 0.707$.  Therefore, our final two points are $p_4 = (0.5, 0.5, 0.707)$ and $p_5 = (0.5, 0.5, -0.707)$.

OPENSCAD is a freely-available program that can be used to generate 3D objects.  For polyhedron (regular or not), the vertices are defined in the code, and then the faces are identified by listing the corners of the face \emph{in clockwise order}.  For our first attempt at creating the octohedron, the OPENSCAD code looks like this:

\begin{verbatim}
polyhedron
    (points = [[0.0, 0.0, 0.0], [1.0, 0.0, 0.0], [1.0, 1.0, 0.0],
        [0.0, 1.0, 0.0], [0.5, 0.5, 0.707], [0.5, 0.5, -0.707],],
        triangles = [[4, 1, 0], [4, 2, 1], [4, 3, 2], [4, 0, 3],
            [5, 0, 1], [5, 1, 2], [5, 2, 3], [5, 3, 0],] );
\end{verbatim}

\indent So, for example, the first face comes from using points $p_4$, $p_1$, and $p_0$, in that order.  The octohedron created by OPENSCAD can be found in Figure 1.

\section{The Octohedron to 3D Print}

While the octohedron created in the previous section is mathematically accurate, it is not usable as created for manufacturing on a 3D printer.  There are two reasons for this.  First, 3D printers are typically programmed using millimeters as its unit of measure, and our first octohedron is just 1 millimeter wide.  This issue is easy to address by scaling the vertices appropriately.  The second issue is more complicated and revolves around the idea that an object is ``printed'' by extruding plastic onto a build plate, starting at $z$=0.  So, objects that we want to print should be sitting on the $xy$-plane.  Translating our initial octohedron so that $p_5$ is touching the $xy$-plane is not the answer, because we want a flat base for the object, and not a sharp point.

We will address these issues by first applying the appropriate rotation to our initial octohedron, and then a scaling factor.  We will rotate about the $x$-axis (which contains $p_0$ and $p_1$ and the rotation will keep those two points fixed).  The rotation angle $\alpha$ will be so that $p_4$ is moved to the $xy$-plane, and with this, all points will satisfy $z \ge 0$.  The rotation matrix that we need is:

\[
R=
  \begin{bmatrix}
    1 & 0 & 0 \\
    0 & \cos \alpha & -\sin \alpha \\
    0 & \sin \alpha & \cos \alpha \\
  \end{bmatrix}
\]

\begin{figure}[h]
\includegraphics[width=\textwidth]{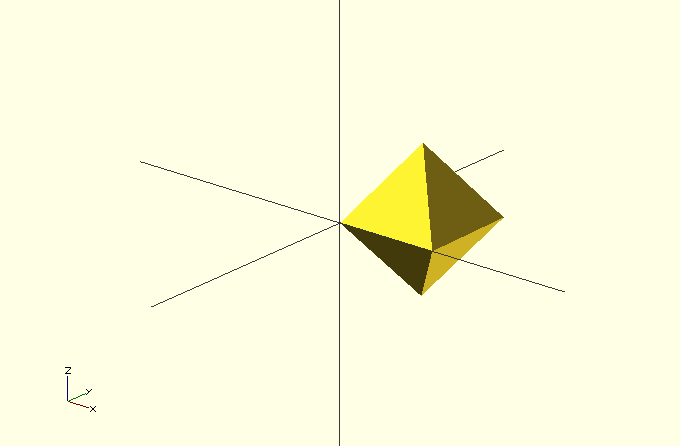}
\caption*{Figure 1: Initial octohedron from OPENSCAD.}
\end{figure}

Multiplying this matrix by $p_0$ or $p_1$ (thought of as vectors) will leave those points fixed.  For $p_4$, we get

\[
Rp_4=
  \begin{bmatrix}
    \frac{1}{2} \\
    \frac{1}{2} \cos \alpha - \frac{\sqrt{2}}{2} \sin \alpha \\
    \frac{1}{2} \sin \alpha + \frac{\sqrt{2}}{2} \cos \alpha \\
  \end{bmatrix}
\]

\noindent In order for this rotated point to sit in the $xy$-plane, we need $\frac{1}{2} \sin \alpha + \frac{\sqrt{2}}{2} \cos \alpha = 0$, which reduces to $\tan \alpha = -\sqrt{2}$.  It follows that $\sin \alpha = \frac{\sqrt{6}}{3} $ and $\cos \alpha = -\frac{\sqrt{3}}{3}$, and $\alpha$ is approximately $-0.955$ radians, or $-54.74$ degrees.

Now that we know $R$, we can apply it to our six points to get the following transformed points (to three decimal places):

\[
\begin{matrix}
\hat{p}_0 = (0.0, 0.0, 0.0) \\
\hat{p}_1 = (1.0, 0.0, 0.0) \\
\hat{p}_2 = (1.0, -0.577, 0.816) \\
\hat{p}_3 = (0.0, -0.577, 0.816) \\
\hat{p}_4 = (0.5, -0.865, 0.0) \\
\hat{p}_5 = (0.5, 0.288, 0.816) \\
\end{matrix}
\]

\noindent Replacing our original points in OPENSCAD with these six points (and using the same triangles) yields another octohedron with the same side lengths, this time sitting on the $xy$-plane with all points satisfying $z \ge 0$.  See Figure 2.

\begin{figure}[h]
\includegraphics[width=\textwidth]{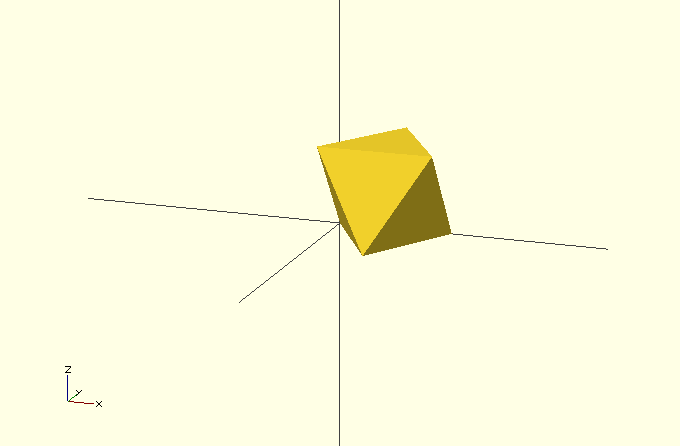}
\caption*{Figure 2: Rotated octohedron from OPENSCAD.}
\end{figure}

Finally, we scale our points in order to manufacture an octohedron of a reasonable size.  Multiplying all values by 20 works out well, and the 3D printed octohedron can be found in Figure 3.

\[
\begin{matrix}
\tilde{p}_0 = (0.0, 0.0, 0.0) \\
\tilde{p}_1 = (20.0, 0.0, 0.0) \\
\tilde{p}_2 = (20.0, -11.54, 16.32) \\
\tilde{p}_3 = (0.0, -11.54, 16.32) \\
\tilde{p}_4 = (10.0, -17.3, 0.0) \\
\tilde{p}_5 = (10.0, 5.76, 16.32) \\
\end{matrix}
\]

\begin{figure}[h]
\center
\includegraphics[width=0.5\textwidth]{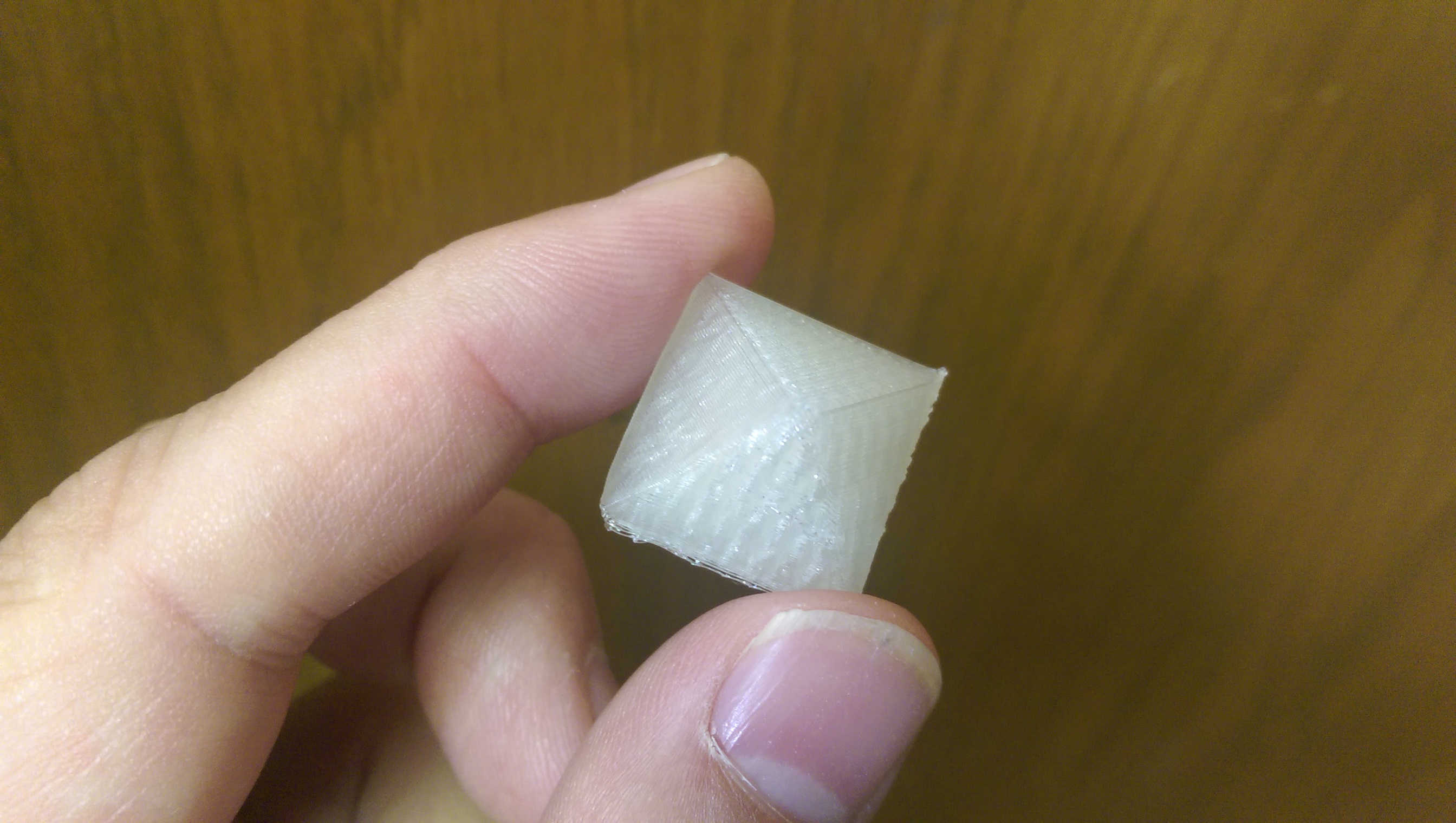}
\includegraphics[width=0.3\textwidth]{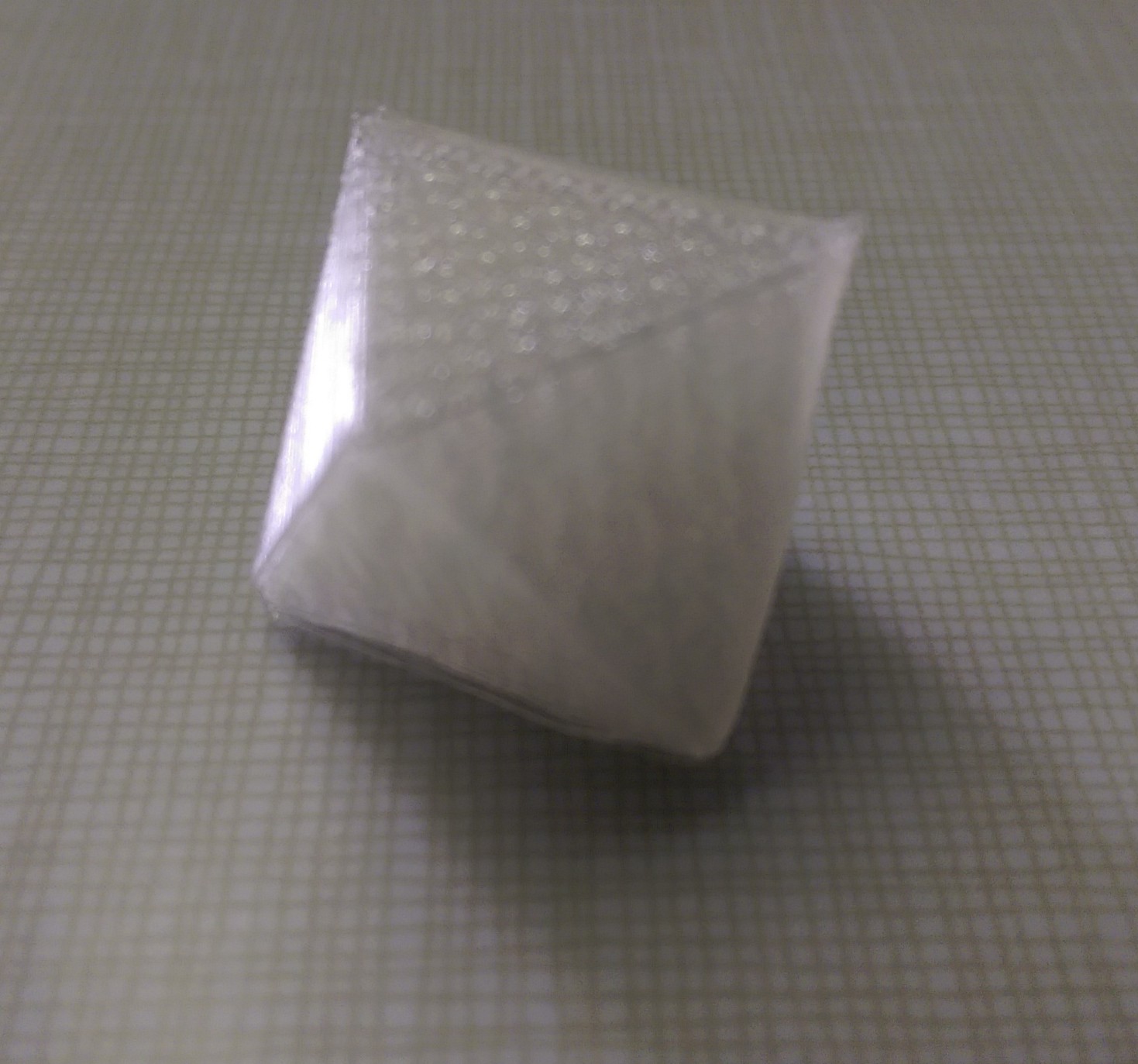}
\caption*{Figure 3: 3D printed octohedron.}
\end{figure}

\section{Appendix: URLs}
\begin{itemize}
\item OPENSCAD home page \url{http://www.openscad.org/}
\item Prof.~Aboufadel's 3D printing page \url{http://sites.google.com/site/aboufadelreu/Profile/3d-printing}
\end{itemize}

\end{document}